\documentclass[11pt]{article}

\usepackage[english]{babel}
\usepackage{amsfonts} 
\usepackage{amsbsy}
\usepackage{amsmath}
\usepackage{amsthm}
\usepackage[dvipdf]{graphicx}
\usepackage{amssymb}
\usepackage{latexsym}
\usepackage{graphicx,color,amsmath,amssymb}
\usepackage{color}

\title{H\"older continuity for optimal multivalued mappings\thanks{
[MS] gratefully acknowledges the University of Toronto
for its generous hospitality during various stages of this work.
[RJM]'s research was supported in part by NSERC grant 217006-08.
\hfill\copyright 2010 by the authors.
}}

\author{R. J. McCann\thanks{Department of Mathematics, University of Toronto, 40 St. George Street, Toronto, Ontario, Canada, M5S 2E4. e-mail: mccann@math.utoronto.ca}\and M. Sosio\thanks{Dipartimento di Matematica ``F. Casorati'',
            Universit\`a degli Studi di Pavia, via Ferrata 1, 27100, Pavia, Italy.  e-mail: maria.sosio@unipv.it}}

\date{}

\begin{document}

\maketitle

\begin{abstract}
Gangbo and McCann showed that optimal transportation between hypersurfaces generally leads to multivalued optimal maps -- bivalent when the target surface is strictly convex. In this paper we quantify H\"older continuity of the bivalent map optimizing average distance squared between arbitrary measures supported on Euclidean spheres.
\end{abstract}



\section{Introduction}

Let $X$, $Y$ be two measure spaces, $\mu$, $\nu$ two probability measures defined on $X$ and $Y$, respectively, and $c$ a measurable map from $X\times Y$ to $[0,+\infty]$. Let us denote with $\Gamma$ the set of all the probability measures on $X\times Y$ that have marginals $\mu$ and $\nu$. The minimization problem
\begin{displaymath}
\mathrm{minimize} \qquad \int_{X\times Y}c(x,y)d\gamma(x,y)\qquad\textrm{among }\gamma\in\Gamma(\mu,\nu)
\end{displaymath}
is known as \emph{Kantorovich's optimal transportation problem}; $c$ is called the \emph{cost function}, and every probability measures in $\Gamma(\mu,\nu)$ is called a \emph{transference plan}. Kantorovich's problem is meant to show how a certain mass $\mu$ distributed on a domain $X$ is transported to another location (described by $\nu$ and $Y$) at a minimal cost (see \cite{villani} for an exhaustive description).

When a measurable map $T$ supports the (unique) optimizer $\gamma\in\Gamma(\mu,\nu)$ for Kantorovich's problem, $T$ is called an \emph{optimal map}.
Several authors treated the regularity of optimal maps when the cost function is the Euclidean squared distance; among them Caffarelli \cite{caff1} \cite{caff2} \cite{caff3} \cite{caff4} \cite{caff5} \cite{caff6}, Delano\"e \cite{del1}, and Urbas \cite{urbas}.
In particular, Caffarelli showed that if the domain $Y$ is convex, $d\mu=f \textrm{\emph{dVol}},d\nu=g\textrm{\emph{dVol}}$, where \emph{dVol} denotes the Lebesgue measure, and the densities $f,1/g$ are bounded, then the optimal map is H\"older continuous.

Recently, in \cite{mtw}, Ma, Trudinger, and Wang extended regularity results to problems with different cost functions by introducing on the cost function a condition (named A3 in their paper), which implies the regularity of the optimal map.
Their work was inspired by Wang's, and Oliker and Waltman's work on the problem of reflector antenna \cite{wang}  \cite{oliker} , which has been shown to be equivalent to optimal transportation of measures on the Euclidean unit sphere with respect to the cost function $-\log|x-y|$ (see for example \cite{wangII} and \cite{glimmoliker}) .

Following Ma, Trudinger, and Wang's paper, Loeper \cite{loeper}, Kim, and McCann \cite{KimMcCann} \cite{KMcCappendices} clarified the role of (A3) in the regularity theory.
More precisely, under the conditions (A0)--(A3), suitable convexity hypotheses on the domains, and the absolute continuity of the Lebesgue measure with respect to $\nu$ \cite{mtw}, Loeper was able to prove the H\"older continuity of the optimal map.
On the other hand, Kim and McCann \cite{KMcCappendices} found a covariant expression of (A3), named (A3s) in their paper, and extended the theory of Loeper to transportation problems set on a pair of smooth manifolds. Very recently, Liu \cite{liu} improved Loeper's H\"older exponent  from $1/(4n-1)$ to $1/(2n-1)$, and showed that the latter is sharp.

Our paper relies on the theory of Loeper, Kim, and McCann to improve the regularity results obtained by Gangbo and McCann \cite{shape} for a transportation problem between boundaries of convex sets. Optimal transportation between boundaries of convex sets does not generally lead to a single-valued map, but rather to multivalued mappings. This is the case of the Kantorovich problem analysed by Gangbo and McCann, who found a bivalent mapping. The novelty of our paper is the quantification of the continuity in this setting of multi-valued mappings.

Let $\Omega$ and $\Lambda$ be two bounded, strongly convex, open sets in $\mathbf{R}^{n+1}$, with Borel probability measures $\mu$ on $\partial\Omega$ and $\nu$ on $\partial\Lambda$. We consider the Monge-Kantorovich problem
\begin{equation}\label{MKprob}
\mathrm{d}^2(\mu,\nu):=\inf_{\gamma\in\Gamma(\mu,\nu)}\int_{\mathbf{R}^{n+1}\times\mathbf{R}^{n+1}}|x-y|^2 d\gamma(x,y),
\end{equation}
where $\Gamma(\mu,\nu)$ denotes the set of all Borel measures on $\mathbf{R}^{n+1}\times\mathbf{R}^{n+1}$ having $\mu$ and $\nu$ as marginals.
Gangbo and McCann \cite{shape} showed that  two maps, named $t^+$ and $t^-$, are required to support the unique optimizer $\gamma\in\Gamma(\mu,\nu)$. This means that the mass at a point $x\in\partial\Omega$ does not always have a unique destination on $\partial\Lambda$, but can be split into two different destinations, $t^+(x)$ and $t^-(x)$. Gangbo and McCann proved that the images of a source $x\in\partial\Omega$ depend  on $x$ in a continuous and continuously invertible way through $t^+$ and $t^-$. In particular, $t^+$ is a homeomorphism between $\partial\Omega$ and $\partial\Lambda$. Moreover, they conjectured H\"older regularity for $t^+$ on $\partial\Omega\setminus S_0$, where
\begin{displaymath}
S_0:=\{x\in\partial\Omega\mid n_\Omega(x)\cdot n_\Lambda(t^+(x))=0\}
\end{displaymath}
represents a part of the `boundary' between the region where the mass splits and the region where it does not. More precisely, if $S_2$ denotes the region where the mass splits (bivalent region), then $S_0$ contains those limit points of $S_2$ at which the split images degenerate to a single image.

The peculiarity of (\ref{MKprob}) is the `hybrid' setting given by the choice of the Euclidean squared distance cost for a transportation problem set on embedded hypersurfaces. One of the difficulties we encountered has been to combine the convexity notion deriving from the Euclidean cost with the dimension and the pseudo-Riemannian structure of the manifolds where the measures are supported. Since the Hausdorff dimension of $\mathrm{spt}\mu$ and $\mathrm{spt}\nu$ is $n$ rather than $n+1$, we are not able to adapt Caffarelli's regularity theory to our problem; (see however \cite{FKM}). Nevertheless Gangbo and McCann's conjecture about H\"older continuity is reinforced by Example 2.4 of Kim-McCann \cite{KimMcCann}: the authors showed that the Euclidean squared distance cost, in the settings of (\ref{MKprob}), satisfies (A0), (A1), (A2) and (A3s) on
\begin{displaymath}
N:=\{(x,y)\in\partial\Omega\times\partial\Lambda \mid n_\Omega(x)\cdot n_\Lambda(y)>0\}.
\end{displaymath}
Despite this comforting result, the regularity of $t^+$ is not immediate. Loeper's theory needs to be adapted to our `hybrid' setting. Moreover the target measure with respect to $t^+$, $\nu_1$, which is the portion of mass `transferred' by $t^+$, does not inherit the hypothesis on $\nu$ of having a positive lower bound on its density with respect to the Lebesgue surface measure. This means there are regions in $\partial\Omega$ where the Lebesgue surface measure is not absolutely continuous with respect to $\nu_1$, so one necessary hypothesis of Loeper's theory is not satisfied. We will treat these regions separately with a different argument.

Our paper is organized as follows. In Section \ref{prel} we report the main result of Gangbo and McCann's paper \cite{shape}; we also discuss the most important statement of this paper and the strategy we are going to adopt to prove it. We will restrict our argument to the case of spherical domains, $\partial\Omega=\partial\Lambda=\mathbf{S}^n$, though we believe that our regularity result can be extended to more general uniformly convex domains. In Section \ref{not} we introduce and clarify some notation. In Section \ref{srq} we comment on some questions related to our problem. In Section \ref{firstpart} we adapt Loeper's theory to our transportation problem, restricting his argument to the regions of $\partial\Omega$ where the necessary hypothesis on the measures holds. The regularity result on the remaining regions is derived in Section \ref{secondpart} .


\section{Preliminaries, strategy, and results}\label{prel}

We recall the following definitions from \cite{shape}. For a smooth convex domain $\Omega$,
          {\em strong} convexity asserts the existence of a positive
          lower bound for all principal curvatures of $\partial \Omega$.
\newtheorem{suitmeas}{Definition}[section]
\begin{suitmeas}\label{suitmeas}
A pair of Borel measures $(\mu,\nu)$ respectively on $(\partial\Omega,\partial\Lambda)$ is said to be {\em suitable} if
\begin{itemize}
\item[(i)] there exists $\epsilon >0$ such that $\mu<\frac{1}{\epsilon}\mathcal{H}^n\lfloor_{\partial\Omega}$ and $\nu>\epsilon\mathcal{H}^n\lfloor_{\partial\Lambda}$, and
\item[(ii)] $\Omega$ is strongly convex.
\end{itemize}
If the above hypotheses are satisfied also when the roles of $\mu\leftrightarrow\nu$ and $\Omega\leftrightarrow\Lambda$ are interchanged, we say the pair $(\mu,\nu)$ is symmetrically suitable.
\end{suitmeas}
Under these hypotheses on the measures, Gangbo and McCann were able to prove the following optimality result.
\newtheorem{opttheo}[suitmeas]{Theorem}
\begin{opttheo}\label{opttheo}
Fix bounded, strictly convex domains $\Omega,\Lambda\in\mathbf{R}^{n+1}$ with suitable measures $\mu$ on $\partial\Omega$ and $\nu$ on $\partial\Lambda$. Then the graphs of a pair of continuous maps $t^+:\partial\Omega\longrightarrow\partial\Lambda$ and $t^-:\bar{S}_2\longrightarrow\bar{T}_2$ contain the support of the unique minimizer $\gamma\in\Gamma(\mu,\nu)$ for \emph{(\ref{MKprob})}:
\begin{equation}
\begin{array}{ccc}
\{(x,t^+(x))\}_{x\in\mathrm{spt}\mu}&\subset\mathrm{spt}\gamma\subset&\{(x,t^+(x))\}_{x\in\partial\Omega}\cup\{(x,t^-(x)))\}_{x\in S_2}\\\quad&\quad&(=\partial\psi\cap(\partial\Omega\times\partial\Lambda)).
\end{array}\nonumber
\end{equation}
Here $\psi$ is Brenier's convex potential for \emph{(\ref{MKprob})}, $S_2=\partial\Omega\setminus\mathrm{dom}\nabla\psi$, $T_2\subset\partial\Lambda$, and $t^+=t^-$ on $\bar{S}_2\setminus S_2$, whereas $t^+(x)-t^-(x)\neq 0$ is an outward normal for $\partial\Omega$ whenever $x\in S_2$.
If $(\mu,\nu)$ are symmetrically suitable, then  $t^+$ and $t^-\lfloor_{\bar{S}_2}$ are homeomorphisms and the minimizer $\gamma$ can be expressed by
\begin{equation}
\gamma=\gamma_1+\gamma_2,\quad \gamma_1=(\mathrm{id}\times t^+)_\sharp\mu_1,\quad\gamma_2=(\mathrm{id}\times t^-)_\sharp\mu_2,\nonumber
\end{equation}
where $\mu_1:=(t^+)^{-1}_\sharp\nu_1$, $\mu_2:=\mu-\mu_1$, and $\nu_1:=\nu\lfloor_{T_2^c}$, with $T_2^c:=\partial\Lambda\setminus t^-(S_2)$ and the identification $t^-(x):=t^+(x)$ for $x\in\partial\Omega\setminus S_2$.
\end{opttheo}

Gangbo and McCann \cite{shape} introduced also a decomposition of the domain $\partial\Omega$, which gives a geometric characterization of $t^+$ and $t^-$.
\newtheorem{charact}[suitmeas]{Proposition}
\begin{charact}\label{charact}
Given $\Omega,\Lambda, (\mu,\nu)$ as in Theorem \ref{opttheo}, we decompose $\partial\Omega=S_0\cup S_1\cup S_2$ into three disjoint sets such that
\begin{eqnarray}
\begin{array}{ccccc}x\in S_0& \iff& t^+(x)=t^-(x)&\textrm{ and }&n_{\Omega}(x)\cdot n_{\Lambda}(t^+(x))=0.\\
x\in S_1& \iff& t^+(x)=t^-(x)&\textrm{ and }&n_{\Omega}(x)\cdot n_{\Lambda}(t^+(x))>0.\\
 x\in S_2& \iff& t^+(x)\neq t^-(x)&\textrm{ and }&n_{\Omega}(x)\cdot n_{\Lambda}(t^+(x))>0,\\ {}&{}&{}&{}&n_{\Omega}(x)\cdot n_{\Lambda}(t^-(x))<0.\end{array}\nonumber\end{eqnarray}
Moreover, the sets $S_0$ and $S_0\cup S_2$ are both closed, whereas $S_1$ is open..
\end{charact}

The partition of $\partial\Omega$ will play an important role in our paper, so it is essential to understand the meaning of these sets. The mass lying on $S_0\cup S_1$ is transferred without splitting to a target set on $\partial\Lambda$ by $t^+$, while the mass lying on $S_2$ splits into two destinations, which are described by $t^+$ and $t^-$. For this reason we will call $S_0$ the \emph{degenerate set}, $S_1$ the \emph{non-degenerate univalent set}, and $S_2$ the \emph{bivalent set}. When the measures $(\mu,\nu)$ are symmetrically suitable, an analogous decomposition of $\partial\Lambda=T_0\cup T_1\cup T_2$ can be introduced (see Definition 3.6 of \cite{shape}). In particular $T_2$ is the bivalent set for the Kantorovich transportation problem (\ref{MKprob}), where $(\Omega,\nu)$ and $(\Lambda,\nu)$ are exchanged, with $(\Lambda,\nu)$ playing the role of the source.

Our aim is to prove that the map $t^+:\partial\Omega\longrightarrow\partial\Lambda$ is H\"older continuous on $S_1$ and $S_2$. One of us has been able to show that $t^+$ satisfies bi-Lipschitz estimates when $n=1$ \cite{2d}, but the argument used cannot be extended to higher dimensions. It heavily relies on the results of Ahmad \cite{ahmad}. Here we are developing a different strategy which works for all $n>1$, when $\partial\Omega,\partial\Lambda=\mathbf{S}^n$. We will proceed in two steps. First we will show that $t^+$ is H\"older continuous on the preimage $\mathbf{S}^n \setminus (t^+)^{-1}(T_2 \cup T_0)$ of the set where \begin{displaymath}\nu_1>\epsilon\mathcal{H}^n\lfloor_{\partial\Lambda},\end{displaymath} where $\epsilon$ is the constant from Definition \ref{suitmeas}. This lower bound on $\nu_1$ allows us to adapt the argument used by Kim and McCann in \cite{KMcCappendices}. On $(t^+)^{-1}(T_2)$, where the lower bound fails, the regularity of $t^+$ will be derived from the H\"older continuity of $t^+$ on $S_2$. In the end we will be able to obtain the following result.

\newtheorem{result}[suitmeas]{Theorem}
\begin{result}[H\"older continuity of multivalued maps outside the degenerate set]\label{result}
If $(\mu,\nu)$ are symmetrically suitable measures on $(\mathbf{S}^n,\mathbf{S}^n)$, $n>1$, then
\begin{displaymath}
t^+\in C^{\frac{1}{4n-1}}_{loc}(S_1)\qquad\textrm{and}\qquad t^+\in C^{\frac{1}{4n-1}}_{loc}(S_2).
\end{displaymath}
\end{result}


\section{Notation}\label{not}

The notation we are going to use is similar to that of \cite{shape} and \cite{KimMcCann}, in particular  we refer to Example 2.4 of \cite{KimMcCann}, with $\partial\Omega=\partial\Lambda=\mathbf{S}^n$, $c:\mathbf{S}^n\times\mathbf{S}^n\rightarrow\mathbf{R}$, $c(x,y)=|x-y|^2$, $N:=\{(x,y)\in\mathbf{S}^n\times\mathbf{S}^n \mid n_{\mathbf{S}^n}(x)\cdot n_{\mathbf{S}^n}(y)>0\}$, and $\hat{N}(x):=\{y\in\mathbf{S}^n \mid (x,y)\in N\}$. We will always use the variable $x$ for points on the source domain $\partial\Omega=\mathbf{S}^n$, and the variable $y$ for points on the target domain $\partial\Lambda=\mathbf{S}^n$.

Let us recall the usual system of local coordinates for the points of $\mathbf{S}^n$
\begin{eqnarray}\begin{array}{cc}
\varphi_i:\mathbf{S}^n\cap \{x\in\mathbf{S}^n|x_i>0\}\rightarrow\mathbf{R}^n, & \varphi_i(x)=(x_1,\ldots,x_{i-1},x_{i+1},\ldots,x_n).\nonumber
\end{array}\end{eqnarray}
Following this example, given $x\in\mathbf{S}^n$ and $y\in\hat{N}(x)$ we can consider a system $\pi_x$ of local coordinates projecting on the hyperplane perpendicular to $x$. In this way both $x$ and $y$ can be represented in local coordinates by means of the same map $\pi_x$
\begin{equation}
x\xrightarrow{\pi_x} X,\qquad y\xrightarrow{\pi_x} Y\nonumber,
\end{equation}
where the capital letters stand for the image of the projection.
To simplify the notation, given a function $F:\mathbf{R}^{n+1}\rightarrow\mathbf{R}$ and a projection $\pi_{x_0}$, whenever $x\in\hat{N}(x_0)$ we will write $F(X)$ to denote $F(\pi_x^{-1}(X))=F(x)$. We will therefore write $\psi(X)$, $c(X,Y)$ instead of $\psi(\pi_{x_0}^{-1}(X))$, $c((\pi_{x_0}^{-1}(X),(\pi_{x_0}^{-1}(Y))$. For example, given $x\in\mathbf{S}^n$ and $y\in N(x)$, by mean of $\pi_x$ we can write
\begin{displaymath}
c(X,Y)=|X-Y|^2+(\sqrt{1-|X|^2}-\sqrt{1-|Y|^2})^2.
\end{displaymath}
In local coordinates, we use the notation $Dc=(\frac{\partial c}{\partial X_1},\ldots, \frac{\partial c}{\partial X_n}$ and $\bar{D}c=(\frac{\partial c}{\partial Y_1},\ldots, \frac{\partial c}{\partial Y_n})$ to denote the partial derivatives. The cross partial derivatives $\bar{D}Dc$ at $(x,y)\in N$ define an unambiguous linear map from vectors at $y$ to covectors at $x$.

Hereafter $d\mathcal{H}^n$ denotes the Hausdorff measure   of dimension $n$, $\mathcal{N}_ \rho(B)$ represents the $\rho$-neighbourhood of a set $B$, and $[Y_0,Y_1]$ indicate the Euclidean segment whose extreme points are $Y_0$ and $Y_1$.

In Section \ref{secondpart} we will use the expression `angle between two vectors $z_1$ and $z_2\in\mathbf{R}^{n+1}$'. The term angle refers to the $\arccos\frac{z_1\cdot z_2}{|z_1||z_2|}$.


\section{Some related questions}\label{srq}

\subsection{Relation between the convex potential $\psi$ and the mappings $t^+,t^-$}

In a Kantorovich problem, when the optimal map exists, it is known to be the gradient of a convex potential. So, what is the relation between the optimal mappings $t^+$, $t^-$, and the convex potential $\psi$? Can we derive any regularity for $\psi$ from Theorem \ref{result}? Gangbo and McCann answered to the first question in Lemma 1.6 of \cite{shape}. Indeed the maps $t^+$ and $t^-$ correspond to the outer and inner trace of $\nabla\psi$, respectively. So we can write the subdifferential of $\psi$ in terms of the optimal mappings: $\partial\psi(x)=[t^+(x),t^-(x)]$ at any boundary point $x\in\partial\Omega$. Moreover, in Corollary 4.4 of \cite{shape}, Gangbo and McCann proved that, when $\Omega$ is bounded strongly convex, $\Lambda$ is bounded strictly convex, and $(\mu,\nu)$ are suitable measures on $\partial\Omega,\partial\Lambda$, then $\psi$ is tangentially differentiable along $\partial\Omega$. This answers the second question. From Theorem \ref{result} it follows immediately that
\begin{displaymath}\psi\in C^{1,\frac{1}{4n-1}}_{loc}\textrm{ on }S_1\subset\mathbf{S}^n,\end{displaymath}i.e. on the non-degenerate univalent set, where $\partial\psi(x)=\{\nabla\psi(x)\}=\{t^+(x)\}$.
Notice that the conclusion of Theorem \ref{result} does not imply $\psi\in C^{1,\frac{1}{4n-1}}_{loc}\textrm{ on }S_2$. Indeed, $\psi$ is not differentiable in the normal direction to the sphere on $S_2$. Nevertheless, choosing the coordinates of Lemma A.1 of \cite{shape}, $\frac{\partial\psi}{\partial x_1}$ exists for $i=2,3,\ldots,n+1$, and  $$\frac{\partial\psi}{\partial x_i}(x)=t^+(x)_i=t^-(x)_i,\textrm{ for }i=2,3,\ldots,n+1,\textrm{ and }x\in\mathbf{S}^n.$$
We conclude that the restriction of $\psi$ to $\mathbf{S}^n$ has a derivative which is H\"older continuous locally on $S_1$ and $S_2$.


\subsection{The regularity of $t^+$ on $S_0$}

Presently we do not have any regularity result, except continuity from \cite{shape}, for $t^+$ on the degenerate set $S_0$. On the contrary, we will see in the statements of Theorem \ref{holder1} and Theorem \ref{holder2} that, on $\mathbf{S}^n\setminus S_0$, close to $S_0$ the H\"older constant of $t^+$ provided by our proof may become very big. Moreover, as noticed in Example 2.4 of \cite{KimMcCann}, the nondegeneracy hypothesis (A2) fails on $S_0$. Therefore, we cannot apply Loeper's argument on $S_0$. On the other hand we believe the set $S_0$ to be small. In dimension $n=1$, with $\Omega$ and $\Lambda$ bounded strictly convex planar domains, Ahmad \cite{ahmad} proved that $S_0$ consists of at most two points.


\subsection{Extending the results to more general domains}

Theorem \ref{result} can be extended to the problem of transporting a measure on a given Euclidean sphere to a measure on any other Euclidean sphere, possibly with a different centre and radius. Indeed, identities (9) and (10) of \cite{shape} indicate how to reduce this more general problem to the case treated in this paper.

Thanks to the results in Example 2.4 of \cite{KimMcCann}, Theorem \ref{holder1} can be extended to the transportation problem where the measures $(\mu,\nu)$ are supported on $(\partial\Omega,\partial\Lambda)$, with $\Omega,\Lambda\subset\mathbf{R}^{n+1}$ bounded convex domains with $C^2$-smooth boundaries. 




\section{$t^+$ is H\"older continuous on $(t^+)^{-1}(T_1)\subset\mathbf{S}^n$}\label{firstpart}

In this section we are going to adapt Kim-McCann's version of Loeper's argument (Appendices B,C and D of \cite{KMcCappendices}) to our mapping $t^+$, which satisfies $(t^+)_\sharp\mu_1=\nu_1$.
Thus, let us recall the regularity conditions (A0),(A1), (A2), and (A3s) from \cite{KimMcCann} \cite{mtw} on a cost function $c:\mathbf{S}^n\times\mathbf{S}^n\rightarrow\mathbf{R}$\newline
\textbf{(A0)(Smoothness)}\emph{ $c\in C^4(N)$, where $N$ has been define in Section \ref{not}.}\newline
\textbf{(A1)(Twist condition)} \emph{A cost $c\in C^1(N)$ is called twisted if for all $x\in \mathbf{S}^n$ the map $y\rightarrow -Dc(x,y)$ from $\hat{N}(x)\subset\partial\Lambda$ to $T_x^*(\mathbf{S}^n)$ is injective.}
\newline
\textbf{(A2)(Non-degeneracy)} \emph{A cost $c\in C^2(N)$ is non-degenerate if for all $(x,y)\in N$ the linear map $\bar{D}Dc:T_y\mathbf{S}^n\rightarrow T_x^*\mathbf{S}^n$ is bijective.}
\newline
\textbf{(A3s)(Strictly regular costs)} \emph{A cost $c\in C^4(N)$ is strictly regular on $N$ if it is non degenerate and for every $(x,y)\in N$
\begin{equation}\label{a3w}
\textrm{sec}_{(x,y)}(p\oplus0)\wedge(0\oplus\bar{p})\ge 0 \textrm{ for all null vectors }p\oplus\bar{p}\in T_{(x,y)}N,
\end{equation} 
and equality in (\ref{a3w})  implies $p=0$ or $\bar{p}=0$ (see \cite{KimMcCann} \cite{mtw} for an explanation of the notation \emph{sec}; we do not need to specify it here, for it will not reappear in the sequel).}

Assuming $(\mu,\nu)$ to be suitable measures on $(\mathbf{S}^n,\mathbf{S}^n)$, in order to apply Kim--McCann's argument we need $\nu_1$ to satisfy
\begin{equation}\label{hypnu}
\textrm{ there exists }\epsilon_1\textrm{ such that } \nu_1>\epsilon_1\mathcal{H}^n\lfloor_{\partial\Lambda}.
\end{equation}
From the definition of $\nu_1$ in Theorem \ref{opttheo} 
we see that $\nu_1$ satisfies (\ref{hypnu}) only outside the bivalent set $T_2\in\partial\Lambda=\mathbf{S}^n$, i.e. outside the set where the image of $t^+$ is bivalent. This is the reason why we can state a regularity result only on a portion of the source domain, $(t^+)^{-1}(T_1)\subset\mathbf{S}^n$.
Hereafter we will assume $n>1$.

\newtheorem{holder1}{Theorem}[section]
\begin{holder1}\label{holder1}
Suppose $(\mu,\nu)$ are symmetrically suitable measures on $(\mathbf{S}^n,\mathbf{S}^n)$ (in particular, from Definition \ref{suitmeas}, there exists $\epsilon>0$ such that $\nu>\epsilon\mathcal{H}^n\lfloor_{\mathbf{S}^n}$). Then $t^+$ is locally H\"older continuous on $(t^+)^{-1}(T_1)$, with H\"older exponent at least $\frac{1}{4n-1}$. Our control on the local H\"older constant depends on $\epsilon$, $n$, and tends to infinity when one approaches the boundary of $N$.
\end{holder1}

\newtheorem{constexplod}[holder1]{Remark}
\begin{constexplod}
Computations that show the explicit dependence of the H\"older constant on the distance of the boundary of $N$ can be found in \cite{phdthesis}.
\end{constexplod}

\newtheorem{bi-convexity}[holder1]{Lemma}
\begin{bi-convexity}\label{bi-c}
The set
$$
N=\{(x,y)\in\mathbf{S}^n\times\mathbf{S}^n \mid n_{\mathbf{S}^n}(x)\cdot n_{\mathbf{S}^n}(y)>0\}
$$
is bi-convex in the sense of Definition 2.5 of \cite{KimMcCann}.
\end{bi-convexity}
{\em Proof:} Fix $x_0\in\mathbf{S}^n$. $\hat{N}(x_0)$ appears convex from $x_0$ if and only if $Dc(x_0,\hat{N}(x_0))$ is convex in $T^*_{x_0}(\mathbf{S}^n)$. Suppose $Dc(x_0,y_0), Dc(x_0,y_1)\in Dc(x,\hat{N}(x))$, where $y_0,y_1\in \hat{N}(x)$. We are going to show that for every $\theta\in (0,1)$
\begin{equation}\label{condconv}
\theta Dc(x_0,y_1)+(1-\theta)Dc(x_0,y_0)\in Dc(x_0,\hat{N}(x)).
\end{equation}
Let's consider a system of local coordinates. Given $x_0\in\mathbf{S}^n$ we project $x_0$ and $y\in\hat{N}(x_0)$ to the hyperplane perpendicular to $\hat{n}_{\Omega}(x_0)$ and containing the origin (notice that this choice of local coordinates is well defined since $\hat{n}_{\Omega}(x_0)\cdot\hat{n}_{\Lambda}(y_k)>0$, when $y_k\in \hat{N}(x_0)$, $k=0,1$)
\begin{equation}
x_0\xrightarrow{\pi_{x_0}} 0, \quad y\xrightarrow{\pi_{x_0}} Y
\end{equation}
so that, in local coordinates,
\begin{eqnarray}
\begin{array}{c}
x_0=(0,1),\quad y=(Y, \sqrt{1-|Y|^2})\\ 
c(X,Y)=|X-Y|^2+(\sqrt{1-|X|^2}-\sqrt{1-|Y|^2})^2.
\end{array}\nonumber
\end{eqnarray}
We easily get 
\begin{equation}
\frac{\partial c}{\partial X_i}(0,Y)=-2Y_i.\nonumber
\end{equation} 
If $v\in T_x(\partial\Omega)$ and $v_i$ are its coordinate with respect to the basis $\frac{\partial}{\partial X_i}$, I can write
\begin{equation}
Dc(v)(x_0,y)=v(c)(x_0,y)=\sum_{i=1}^n v_i\frac{\partial c}{\partial X_i}(0,Y).\nonumber
\end{equation}
Hence we can compute
\begin{eqnarray}\label{comp}
\begin{array}{c}
\theta Dc(v)(x_0,y_1)+(1-\theta)Dc(v)(x_0,y_0)
=\theta v(c)(x_0,y_1)+(1-\theta)v(c)(x_0,y_0)\\
=\sum_{i=1}^n\left[\theta v_i\frac{\partial c}{\partial X_i}(0,Y_1)+(1-\theta)v_i\frac{\partial c}{\partial X_i}(0,Y_0)\right]\\
=\sum_{i=1}^n 2v_i\left[\theta\left(-Y_{1,i}\right)+(1-\theta)\left(-Y_{0,i}\right) \right]\\
=\sum_{i=1}^n -2v_i(\theta Y_{1,i}+(1-\theta)Y_{0,i}).
\end{array}
\end{eqnarray}
Therefore, for all $\theta\in(0,1)$ 
\begin{equation}
\theta Dc(x_0,y_1)+(1-\theta)Dc(x_0,y_0)= Dc(x_0,\pi_{x_0}^{-1}(\theta Y_1+(1-\theta)Y_0)\in Dc(x_0,\hat{N}(x_0)). \nonumber
\end{equation} 
Since $x_0$ is an arbitrary point of $\mathbf{S}^n$, we conclude that $N$ is vertically convex.
By a similar argument, it is easy to show that $N$ is also horizontally convex. We conclude that $N$ is bi-convex.
$\qquad\square$
\newline

{\em Proof of Theorem \ref{holder1}: }Fix $(x,y)= (x,t^+(x)) \in N$, with $t^+(x) \in T_1$. Since $T_1$ is open,  and $t^+$ is continuous,
we can choose $R$ and then $r$ small enough that $B_r(y) \subset t^+(B_R(x)) \subset T_1$; as asserted by Trudinger and Wang in \cite{TW}, since $N$ is bi-convex, taking  $R$ and $r$ even smaller, $P=B_R(x) \times B_r(y) \subset N$ is bi-convex (alternatively, we could show directly that $P$ is bi-convex, by means of the same argument used for $N$ in Lemma \ref{bi-c}. 
We replace $\nu_1$ with its restriction $\nu'_1$ to $B_r(y)$ and
we denote $\mu'_1 = s^+_\#(\nu'_1)$.
Taking $R$ and then $r$ even smaller than before gives us local coordinates over both domains simultaneously (for example through the chart $\pi_x$).
Let $X=\pi_x(x),Y=\pi_x(y)$, and $P' = \pi_x(B_R(x)) \times \pi_x(B_r(y))$. Since $P$ is bi-convex and the notion of bi-convexity is coordinate invariant (as manifest from Definition 2.5 of \cite{KimMcCann}), $P'$ is bi-convex with respect to the cost
\begin{equation}\label{costlc}
c(X,Y)=|X-Y|^2+(\sqrt{1-|X|^2}-\sqrt{1-|Y|^2})^2.
\end{equation}
Kim and McCann showed that the cost in the original coordinates satisfies condition (A2) and (A3s) (see Example 2.4 of \cite{KimMcCann}), and that the quantities in these conditions have an intrinsic meaning independent of coordinates, since they are geometric quantities (i.e. pseudo-Riemannian curvatures in the case of (A3s) and non-degeneracy of the metric in the case (A2)). This implies that also the cost (\ref{costlc}) satisfies (A2) and (A3s). Only the constant $C'_0$ of (A3s) will depend on the coordinates. Since we know that the equation $D_X c(X,Y) = D\psi(X)$ has at most two solutions, $Y^+=t^+(X)$ and $Y^-=t^-(X)$ and only $Y^+$ lies in $P'$, the cost satisfies (A1) on $P'$.

At this point we can apply Theorem D.1. of \cite{KMcCappendices} to the cost (\ref{costlc}) on $P'$, with probability measures $\mu_1^x$ and $\nu_I^x$, on $\pi_x(B_R(x))$ and $\pi_x(B_r(y))$ respectively, defined by
\begin{equation}
\mu_1^x:=(\pi_x)_\sharp \mu_1',\qquad \nu_1^x:=(\pi_x)_\sharp \nu_1'.\nonumber
\end{equation}
The source $\mu_1^x$ is supported (and bounded above) in $\pi_x(B_R(x))$ and target $\nu_1^x$ supported (and bounded below) in $\pi_x(B_r(y))$,
We deduce the existence of a locally H\"older continuous optimal
map pushing $\mu'_1$ forward to $\nu'_1$.  By the uniqueness of
optimal transport, this map must coincide
$\mu'_1$-a.e. with $t^+$.  Since both maps are continuous they
agree on the (closed) support of $\mu'_1$.
Since $\mathrm{spt} \mu'_1$ contains a small ball around $x$, this shows $t^+$ is locally H\"older at $x$.
$\qquad\square$


\section{$t^+$ is locally H\"older continuous where its image is bivalent}\label{secondpart}

The previous section established local H\"older continuity
for the outer map $t^+=(s^+)^{-1}$ on the source domain $s^{+}(T_1) \subset \mathbf{S}^n$,
but not on $s^+(T_0 \cup T_2) = S_0 \cup s^+(T_2)$.
Our strategy for extending this estimate to $s^+(T_2)$ is described
at the end of this paragraph. First note, however, that Gangbo and McCann's
\emph{Sole Supplier Lemma}, 2.5 of \cite{shape}, implies the outer
image of the bivalent source is contained in the univalent target
$t^+(S_2) \subset T_1$, and similarly 
$s^+(T_2) \subset S_1$.  Since $s^+: \mathbf{S}^n \longrightarrow \mathbf{S}^n$
is a homeomorphism,  from $S_1 \cup S_2 = s^+(T_1) \cup s^+(T_2)$, it follows
that the bivalent source $S_2 \subset s^+(T_1)$ belongs to the domain where
H\"older continuity of $t^+$ has already been shown. On this bivalent set $S_2$,
the inner map $t^-$ is related to the outer map $t^+(x) = t^-(x) + \lambda (x) x$
by the geometry of the target.  In Proposition \ref{t-holder}, this relation will be used
to deduce (i) H\"older continuity of $t^-$
from that of $t^+$. 
This quantifies injectivity (ii) of the inverse map
$s^- = (t^-)^{-1}$ (through a bi-H\"older estimate in Proposition \ref{s-holderbelow}),
whose relation to the outer map $s^+(y) = s^-(y) + \omega(y)y$
is then used in Proposition \ref{s+holderbelow}
to quantify injectivity (iii) of $s^+ = (t^+)^{-1}$ on the bivalent target
$T_2 = t^-(S_2)$.
This yields the desired local H\"older continuity of $t^+$ on the
source set $s^+(T_2)$ mentioned at the outset.

Let us recall the geometric characterization of $t^+$ and $t^-$ from Definition \ref{charact}. Remembering that, on $\mathbf{S}^n$, $n_{\mathbf{S}^n}(x)=x$, we have
\begin{itemize}
\item If $x\in S_0$ then
$x\cdot t^+(x)=0$.
\item If $x\in S_1$ then
$x\cdot t^+(x)>0$.
\item If $x\in S_2$ then
$x\cdot t^+(x)>0$ and $x\cdot t^-(x)<0$.
\end{itemize}
We are going to introduce a geometric approach, based on the previous characterization, which allows us to prove the following theorem. Hereafter we will assume $n>1$.

\newtheorem{holder2}{Theorem}[section]
\begin{holder2}\label{holder2}
If $(\mu,\nu)$ are symmetrically suitable measures on $(\mathbf{S}^n,\mathbf{S}^n)$, then $t^+$ is locally H\"older continuous on $(t^+)^{-1}(T_2)$.
\end{holder2}

From Lemma 1.6 of \cite{shape} we know that $t^+$ and $t^-$ are related by
\begin{equation}
\forall x\in S_2\subset\mathbf{S}^n\qquad t^+(x)-t^-(x)=\lambda(x)x,\nonumber
\end{equation}
where $\lambda$ is a continuous positive function on $S_2$. Given $x_0,x_1$ in $S_2$ we then have
\begin{equation}\label{t+-}
|t^-(x_1)-t^-(x_0)|\le|t^+(x_1)-t^+(x_0)-\lambda(x_1)x_1+\lambda(x_0)x_0|.
\end{equation}
We would like to exploit the regularity of $t^+$ on $S_2 \subset (t^+)^{-1}(T_1)$, proved in the previous section, to prove that also $t^-$ is H\"older continuous on $S_2$. For this purpose we also need to estimate the term $\lambda(x_1)x_1+\lambda(x_0)x_0$. This will be done applying the Mean Value Theorem to a suitable function and utilizing the geometric properties of the target.

\newtheorem{t-holder}[holder2]{Proposition}
\begin{t-holder}[H\"older continuity of $t^-$]\label{t-holder}

If $t^+\in C^{\alpha}_{loc}(S_2)$ then $t^-\in C^{\alpha}_{loc}(S_2)$.
Let $U\subset S_2$ and $0<k_U:=\min\{-x\cdot t^+(x)\mid x\in U\}$. If $C_U^+$ bounds the H\"older constant for $t^+$ on $U$, then \begin{displaymath}C^-_U:=\left(1+\frac{1}{k_U}\right)(C^+_U+2)\end{displaymath} is the H\"older constant for $t^-$ on $U$.

\end{t-holder}
{\em Proof: }The function $h(y):= d(y,\mathbf{S}^n)=1-|y|$ is differentiable on $\Lambda=B_1(0)$ except at $y=0$.
Notice that $h(t^-(x))=h(t^+(x)-\lambda(x)x)=h(t^+(x))=0$ whenever $x\in S_2$.
Consider a neighbourhood $U\subset S_2$ and the corresponding $k_U, C^+_U$ from the statement of Proposition \ref{t-holder}.
Let $x_0,x_1\in U$, $|x_1-x_0|<2$ (we need $\nabla h$ to be well defined on the line segment between $t^-(x_0)$ and $t^-(x_1)$, i.e. $0\notin[t^-(x_0,t^-(x_1)]$).
Applying the Mean Value Theorem, we get
\begin{eqnarray}
\lefteqn{0=h(t^+(x_1)-\lambda(x_1)x_1)-h(t^+(x_0)-\lambda(x_0)x_0)}\nonumber\\
&=& \nabla h(u)\cdot(t^+(x_1)-t^+(x_0)-\lambda(x_1)x_1+\lambda(x_0)x_0)\nonumber,
\end{eqnarray}
for some $u$ on the line segment between $t^-(x_0)$ and $t^-(x_1)$. It follows
\begin{equation}\label{mvt}
(\lambda(x_1)x_1-\lambda(x_0)x_0)\cdot\nabla h(u)=(t^+(x_1)-t^+(x_0))\cdot\nabla h(u).
\end{equation}
We can rewrite (\ref{mvt}) as
\begin{eqnarray}
(t^+(x_1)-t^+(x_0))\cdot\nabla h(u)+\lambda(x_0)(x_0-x_1)\cdot\nabla h(u)\nonumber\\
=(\lambda(x_1)-\lambda(x_0))x_1\cdot\nabla h(u);\nonumber
\end{eqnarray}
then, using $|\nabla h(u)|=1$,
\begin{eqnarray}\label{tiq}
\lefteqn{|\lambda(x_1)-\lambda(x_0)||x_1\cdot\nabla h(u)|}\nonumber\\
&\le& |(t^+(x_1)-t^+(x_0))|+\lambda(x_0)|x_0-x_1|.
\end{eqnarray}

We now state a claim, whose demonstration  is postponed to the end of this proof.
\newtheorem{normali}[holder2]{Lemma}
\begin{normali}\label{normali}
Under the hypotheses of Proposition \ref{t-holder}, fix $\epsilon\in(0,1)$, such that $\epsilon^2<\frac{k_U}{2}$.
Since $t^-$ is uniformly continuous on $\bar{S_2}$, there exists $\delta_\epsilon$, depending on the data through $\psi$, such that
\begin{displaymath}
|x_1-x_0|<\delta_\epsilon\Rightarrow|t^-(x_1)-t^-(x_0)|<\epsilon.
\end{displaymath}
Then,
taking $x_0,x_1$ such that $|x_1-x_0|<\delta_\epsilon$, we have
\begin{displaymath}x_i\cdot\nabla h(u)>\frac{k_U}{2}>0\quad\quad\textrm{ for }i=1,2.\end{displaymath}
\end{normali}

Recalling that $\lambda(x)\le2$, since $\partial\Omega=\mathbf{S}^n$, by means of Lemma \ref{normali} we simplify (\ref{tiq}) to
\begin{eqnarray}\label{stimalambda}
\lefteqn{|\lambda(x_1)-\lambda(x_0)|}\nonumber\\&\le&\frac{2}{k_U}\left[ |t^+(x_1)-t^+(x_0)|+\lambda(x_0)|x_0-x_1|\right]\nonumber\\
&\le& \frac{2}{k_U} [|t^+(x_1)-t^+(x_0)|+2|x_1-x_0|].
\end{eqnarray}
Therefore, by (\ref{t+-}) and (\ref{stimalambda}),
\begin{eqnarray}\label{t-t+}
\lefteqn{|t^-(x_1)-t^-(x_0)|}\nonumber\\&\le&|t^+(x_1)-t^+(x_0)|+\lambda(x_1)|x_1-x_0|+|\lambda(x_1)-\lambda(x_0)|\nonumber\\
&\le& |t^+(x_1)-t^+(x_0)|+2|x_1-x_0|+|\lambda(x_1)-\lambda(x_0)|\nonumber\\
&\le&\left(1+\frac{2}{k_U}\right)|t^+(x_1)-t^+(x_0)|+2 \left(1+\frac{2}{k_U}\right)|x_1-x_0|.
\end{eqnarray}
Combining (\ref{t-t+}) and $t^+\in C^{\alpha}(U)$, we conclude \begin{eqnarray}\label{holderlip}
\lefteqn{|t^-(x_1)-t^-(x_0)|}\\
&\le& C^+_U\left(1+\frac{2}{k_U}\right)|x_1-x_0|^\alpha+2 \left(1+\frac{2}{k_U}\right)|x_1-x_0|\nonumber,
\end{eqnarray}
i.e. $t^-$ is H\"older continuous on $S_2$ whenever $|x_1-x_0|<\delta_\epsilon$, with $\epsilon^2<\frac{k_U}{2}$.
We can take $\delta_\epsilon<1$, so that (\ref{holderlip}) implies
\begin{eqnarray}
|t^-(x_1)-t^-(x_0)|
&\le& \left(1+\frac{2}{k_U}\right)\left[C^+_U+2\right]|x_1-x_0|^\alpha \nonumber\\&=&C^-_U|x_1-x_0|^\alpha\nonumber.\qquad\square
\end{eqnarray}
%

{\em Proof of Lemma \ref{normali}: }
Let $z_i=t^-(x_i),i=1,2$.
Notice that $\nabla h(u)=-\frac{u}{|u|}$. We have $u=sz_1+(1-s)z_0$ for some $s\in(0,1)$. Hence, there exists $\xi\in(0,\epsilon)$ such that
\begin{eqnarray}
\lefteqn{x_1\cdot u<-k_Us+(1-s)x_1\cdot z_0}\nonumber\\
&=&-k_Us+(1-s)x_1\cdot(z_1+\xi(z_0-z_1))\nonumber\\
&<&-k_U+(1-s)\xi\epsilon<-k_U+\epsilon^2.\nonumber
\end{eqnarray}
Using a similar argument for $x_0\cdot u$, we conclude that if $\epsilon^2<\frac{k_U}{2}$ then  $x_i\cdot\nabla h(u)>\frac{k_U}{2|u|}>\frac{k_U}{2}>0$, for $i=1,2$.$\qquad\square$


\newtheorem{viceversa}[holder2]{Remark}
\begin{viceversa}
Proposition \ref{t-holder} admits a converse, i.e. if $t^-\in C^{\alpha}_{loc}(S_2)$ then $t^+\in C^{\alpha}_{loc}(S_2)$. This can be proved with minor changes in the preceding argument.
\end{viceversa}


\newtheorem{t-holderindeed}[holder2]{Remark}
\begin{t-holderindeed}
By means of Theorem \ref{holder1} and Proposition \ref{t-holder}, $t^-$ is indeed locally H\"older continuous on $S_2$ with exponent $\frac{1}{4n-1}$.
\end{t-holderindeed}


The injectivity (ii) of the inverse map
$s^- = (t^-)^{-1}$ on $T_2$, is an immediate consequence of the local H\"older continuity of $t^-$ on $S_2$, and it has been included in the following proposition.

\newtheorem{s-holderbelow}[holder2]{Proposition}
\begin{s-holderbelow}[Quantifying injectivity of $s^-$]\label{s-holderbelow}
Let $V\subset T_2$. Under the hypotheses of Theorem \ref{holder2}
$s^-:=(t^-)^{-1}$ satisfies\begin{displaymath}\forall y_0,y_1\in V \textrm{ sufficiently close}, | s^-(y_1)-s^-(y_0)|\ge \hat{C}^-_V|y_1-y_0|^{4n-1},\end{displaymath}
where
\begin{displaymath}
\hat{C}^-_V=(C^-_U)^{-1},
\end{displaymath}
with $U=s^-(V)$ and $0<k_V:=\min\{-y\cdot s^-(y)\mid y\in V\}$.
\end{s-holderbelow}
{\em Proof: }
Since $s^-:=(t^-)^{-1}$ is uniformly continuous on $\bar{T}_2$, given $\delta_{\epsilon}>0$ there exists $\gamma_{\delta_{\epsilon}}>0$ such that, if $|y_1-y_0|<\gamma_{\delta_{\epsilon}}$, then $|s^-(y_1)-s^-(y_0)|<\delta_\epsilon$. Supposing $|y_1-y_0|<\gamma_{\delta_{\epsilon}}$, we can apply Proposition \ref{t-holder} to $x_1=s^-(y_1), x_0=s^-(x_0)$ to get
\begin{equation}
|s^-(y_1)-s^-(y_0)|
\ge\frac{1}{C^-_U}|y_1-y_0|^{4n-1}.\nonumber\qquad\square
\end{equation}


We now state an elementary Lemma about vectors in $\mathbf{R}^n$.
\newtheorem{vectors}[holder2]{Lemma}
\begin{vectors}\label{vectors}
Let $u,v\in \mathbf{R}^n$. Suppose the angle between $u$ and $v$ is less than $\frac{\pi}{2}+\alpha$, with $\alpha\in\left[0,\frac{\pi}{2}\right)$. Then $|u+v|\ge|u|\cos\alpha$.
\end{vectors}
{ \em Proof: } Let $\theta_{u,v}$ denote the angle between $u$ and $v$.
Keeping $|u|$ and $|v|$ fixed, $|u+v|$ can be seen as a function of $\theta_{u,v}$ by mean of
\begin{displaymath}
|u+v|^2\left(\theta_{u,v}\right)=|u|^2+|v|^2+2|u||v|\cos\theta_{u,v},
\end{displaymath}
When  $\theta_{u,v}\in\left[0,\frac{\pi}{2}+\alpha\right]$, the function $|u+v|\left(\theta_{u,v}\right)$ reaches its minimum at $\theta_{u,v}=\frac{\pi}{2}+\alpha$. To our purpose we can take $\theta_{u,v}=\frac{\pi}{2}+\alpha$.
For simplicity we assume $v$ parallel to $e_1\in \mathbf{R}^n$. Let us consider the projection $p$ on the hyperplane perpendicular to $e_1$ and containing the origin. Then $p(u+v)=p(u)=|u|\cos\alpha$. Since $|p(u+v
)|\le|u+v|$, we have the thesis. $\qquad\square$

This Lemma turns out to be the key to the proof of step (iii).
Under the hypothesis of symmetrically suitable measures, the optimal transportation problem we are studying is symmetric, hence every result that holds for $t^+$ on $\mathbf{S}^n$ implies an analogous result for $s^+$ on $\mathbf{S}^n$. In particular, from Lemma 1.6 of \cite{shape}, for every $y\in T_2$ we can write
\begin{equation}\label{s+-}
s^+(y)-s^-(y)=\omega(y)y,
\end{equation}
where $\omega$ is a nonnegative function on $T_2$.
Hence
\begin{equation}
|s^+(y_1)-s^+(y_0)|=|s^-(y_1)-s^-(y_0)+\omega(y_1)y_1-\omega(y_0)y_0|.\nonumber
\end{equation}
If we were allowed to apply Lemma \ref{vectors} to the right hand side of the previous equality, with $u=s^-(y_1)-s^-(y_0)$ and $v=\omega(y_1)y_1-\omega(y_0)y_0$, we would then be able to exploit the regularity of $s^-$ to prove step (iii). Therefore,
we need to understand the behaviour of the angle between $s^-(y_1)-s^-(y_0)$ and $\omega(y_1)y_1-\omega(y_0)y_0$, when $y_0$ gets close to $y_1$.
From the monotonicity of $\partial\psi$ we have
\begin{equation}
(s^-(y_1)-s^-(y_0))\cdot(y_1-y_0)\ge0\qquad\forall y_1,y_0\in T_2,\nonumber
\end{equation}
which says that the angle between $s^-(y_1)-s^-(y_0)$ and $y_1-y_0$ is in $\left[0,\frac{\pi}{2}\right]$. If we can show that the angle between $y_1-y_0$ and $\omega(y_1)y_1-\omega(y_0)y_0$ is in $[0,\alpha]$, for a certain $\alpha\in\left[0,\frac{\pi}{2}\right)$, then we can apply Lemma \ref{vectors} to get the desired estimate on $|s^+(y_1)-s^+(y_0)|$.

\newtheorem{defbeta}[holder2]{Lemma}
\begin{defbeta}\label{defbeta}
Given $y_0,y_1\in T_2$ we denote with $\beta(y_0,y_1)$ the angle between  $y_1-y_0$ and $\omega(y_1)y_1-\omega(y_0)y_0$. If the angle between $y_0$ and $y_1$ is equal to $\gamma$ then \begin{equation}\label{boundbeta}\beta(y_0,y_1)\in\left[0,\frac{\pi-\gamma}{2}\right).\end{equation}
\end{defbeta}

{\em Proof: } The angle between $y_1$ and $-y_0$ is equal to $\pi-\gamma$, while the angle between $y_1$ (or $-y_0$) and $y_1-y_0$ is $\frac{\pi-\gamma}{2}$. Since $\omega(y_0),\omega(y_1)>0$, $\beta(y_0,y_1)\in\left[0,\frac{\pi-\gamma}{2}\right)$.$\qquad\square$

\newtheorem{anglebeta}[holder2]{Lemma}
\begin{anglebeta}[Dichotomy]\label{anglebeta}
Fix $y_1\in T_2$. For every integer $m>1$ define
\begin{displaymath}
\Theta_m(y_1):=\left\{y\in T_2\mid\beta(y,y_1)\in\left[\frac{\pi}{2}-\frac{1}{m},\frac{\pi}{2}\right]\right\}.
\end{displaymath}
Unless $\Theta_m(y_1)$ is empty for $m$ sufficiently large,
there exist $M>0$ and $K>0$ such that
\begin{equation}\label{hypsup}
|s^+(y_1)-s^+(y)|\ge K|y_1-y|,\qquad\forall y\in \Theta_m(y_1),\textrm{ with }m>M.
\end{equation}
\end{anglebeta}

\noindent
{\em Proof:}
We are interested in the sets $\Theta_m(y_1)$ for $m$ large, so hereafter we assume $m>50$.
Define
\begin{equation}
0< \varpi_m:=\inf \left\{\omega(y)>0\mid y\in\Theta_m(y_1)\right\}\nonumber\end{equation}
and note $\varpi_m \le \varpi_{m+1}$ since $\Theta_m(y_1) \supset \Theta_{m+1}(y_1)$.
By elementary computations, we have
\begin{eqnarray}\label{cos}
\lefteqn{|\omega(y_1)y_1-\omega(y)y|\cos\beta(y,y_1)=\frac{(\omega(y_1)y_1-\omega(y)y)\cdot(y_1-y)}{|y_1-y|}}\nonumber\\
&=&\frac{\omega(y_1)y_1\cdot(y_1-y)-\omega(y)y\cdot(y_1-y)}{|y_1-y|}\qquad\qquad\nonumber\\
&\ge&\frac{\varpi_m y_1\cdot(y_1-y)-\omega(y)y\cdot(y_1-y)}{|y_1-y|}\nonumber\\
&=&\frac{\varpi_m |y_1-y|^2+(\varpi_m-\omega(y))y\cdot(y_1-y)}{|y_1-y|}\nonumber\\
&\ge&\varpi_m|y_1-y|\qquad\qquad \qquad\forall y\in\Theta_m(y_1),
\end{eqnarray}
where we used the definition of $\varpi_m$ and the trivial inequality $y\cdot y_1\le1$ to show that the term $(\varpi_m-\omega(y))y\cdot(y_1-y)$ is non-negative.
Consider now the two vectors $\omega(y_1)y_1-\omega(y)y$ and $(\omega(y_1)-\omega(y))y_1$, with $y\in T_2$. Their difference is parallel to $y_1-y$, so they have the same projection on any hyperplane perpendicular to $y_1-y$. This projection has length $|\omega(y_1)y_1-\omega(y)y|\sin\beta(y,y_1)$. Therefore
\begin{equation}\label{sin}
|\omega(y_1)y_1-\omega(y)y|\sin\beta(y,y_1)\le|\omega(y_1)-\omega(y)|\qquad\forall y\in T_2.
\end{equation}
Putting together (\ref{cos}) and (\ref{sin}), we obtain an estimate for $\tan\left(\frac{\pi}{2}-\frac{1}{m}\right)$
\begin{eqnarray}
\tan\left(\frac{\pi}{2}-\frac{1}{m}\right)\le\tan{\beta(y,y_1)}\le\frac{|\omega(y_1)-\omega(y)|}{\varpi_m|y_1-y|}\qquad\forall y\in\Theta_m(y_1).\nonumber
\end{eqnarray}

As $m\rightarrow+\infty$, $\tan\left(\frac{\pi}{2}-\frac{1}{m}\right)\rightarrow +\infty$; then for every $N>0$ there exists $m_N>50$ such that
\begin{equation}\label{contrad}
|\omega(y_1)-\omega(y)|>N\varpi_m|y_1-y|,\qquad\forall y\in \Theta_m(y_1), m>m_N.
\end{equation}

From (\ref{s+-}) we have, for every $y\in T_2$,
\begin{eqnarray}
s^+(y_1)-s^+(y)-\omega(y)(y_1-y)=
s^-(y_1)-s^-(y)+(\omega(y_1)-\omega(y))y_1.\nonumber
\end{eqnarray}
We define
\begin{eqnarray}\label{defA}
A&:=&|s^+(y_1)-s^+(y)|+|\omega(y)(y_1-y)|\nonumber\\
&\ge&|s^-(y_1)-s^-(y)+(\omega(y_1)-\omega(y))y_1|,\qquad y\in T_2.
\end{eqnarray}
Using $|v-u|\ge|v|-|u|\quad \forall u,v\in\mathbf{R}^{n+1}$, we get two different estimates for $A$
\begin{equation}\label{s1}
A\ge |s^-(y_1)-s^-(y)|-|\omega(y_1)-\omega(y)|,\end{equation}
\begin{equation}\label{s2}
A\ge |\omega(y_1)-\omega(y)|-|s^-(y_1)-s^-(y)|.\end{equation}
By the symmetry of the problem, using (\ref{stimalambda}), we have
\begin{equation}
|\omega(y_1)-\omega(y)|\le\frac{2}{k'_m} \left[|s^+(y_1)-s^+(y)|+2|y_1-y|\right]\qquad\forall y\in\Theta_m(y_1),\nonumber
\end{equation}
where $0<k'_m:=\inf \left\{-y\cdot s^-(y)\mid y\in\Theta_m(y_1)\right\} \le k'_{m+1}$.
From (\ref{s1}) it follows
\begin{equation}\label{s1f}
A\ge |s^-(y_1)-s^-(y)|-\frac{2}{k'_m} \left[|s^+(y_1)-s^+(y)|+2|y_1-y|\right].
\end{equation}
On the other hand, combining (\ref{contrad}) and (\ref{s2})
\begin{equation}\label{s2f}
A\ge N\varpi_m|y_1-y|-|s^-(y_1)-s^-(y)|, \qquad\forall y\in\Theta_m(y_1), m>{m}_N.
\end{equation}
We can sum (\ref{s1f}) and (\ref{s2f}) to get
\begin{equation}
2A\ge N\varpi_m|y_1-y|-\frac{2}{k'_m} \left[|s^+(y_1)-s^+(y)|+2|y_1-y|\right].\nonumber
\end{equation}
From the definition (\ref{defA}) of $A$, this becomes
\begin{equation}
2\left(1+\frac{1}{k'_m}\right)|s^+(y_1)-s^+(y)|\ge\left(N\varpi_m-\frac{4}{k'_m}-2\omega(y)\right)|y_1-y|,\nonumber
\end{equation}
for every $y\in\Theta_m(y_1), m>m_N$.
Since neither $\varpi_m$ nor $k'_m$ is decreasing as a function of $m$,
taking $N$ large enough ensures
$N>\left(\frac{4}{k'_{m_N}}+4\right)\frac{1}{\varpi_{m_N}}$
to yield a positive constant
\begin{displaymath}
K=\frac{N\varpi_{m_N}-4(\frac{1}{k'_{m_N}}+1)}{2\left(1+\frac{1}{k'_{m_N}}\right)},\end{displaymath}
such that
\begin{displaymath}
|s^+(y_1)-s^+(y)|\ge K|y_1-y|,\qquad\forall y\in\Theta_m(y_1),m>m_N.
\end{displaymath}
To conclude we take $M=m_N$.
$\qquad\square$

The injectivity (iii) of $s^+ = (t^+)^{-1}$ on the bivalent target
$T_2 = t^-(S_2)$.
 follows from Lemma \ref{vectors} and Lemma \ref{anglebeta}.

\newtheorem{s+holderbelow}[holder2]{Proposition}
\begin{s+holderbelow}[Quantifying injectivity of $s^+$ on the bivalent target]\label{s+holderbelow}Let $y_1\in V\subset T_2$. Under the hypotheses of Theorem \ref{holder2}, there exists $\hat{C}^+_V>0$, depending on $\hat{C}^-_V$, $k_V$ (from Proposition \ref{s-holderbelow}), and $\bar{\theta}(y_1)$ (from Lemma \ref{anglebeta}), such that, when $y_0$ is sufficiently close to $y_1$,
\begin{displaymath} | s^+(y_1)-s^+(y_0)|\ge \hat{C}^+_V
|y_1-y_0|^{4n-1}.\end{displaymath}
\end{s+holderbelow}

{\em Proof: }
When $y_0\in\Theta(y_1,\theta)$, with $\theta>\bar{\theta}(y_1)$ we apply Lemma \ref{anglebeta} and we are done. Otherwise
the angle between $s^-(y_1)-s^-(y_0)$ and $\omega(y_1)y_1-\omega(y_0)y_0$ is smaller than $\frac{\pi}{2}+\bar{\theta}(y_1)$. Applying Lemma \ref{vectors}, we obtain
\begin{eqnarray}
\lefteqn{|s^+(y_1)-s^+(y_0)|}\nonumber\\&=&|s^-(y_1)-s^-(y_0)+\omega(y_1)y_1-\omega(y_0)y_0|\nonumber\\
&\ge&|s^-(y_1)-s^-(y_0)|\cos{\bar{\theta}(y_1)}.\nonumber
\end{eqnarray}
Taking $y_0,y_1$ sufficiently close ($|y_1-y_0|<\gamma_{\delta_\epsilon}$, from the proof of Proposition \ref{s-holderbelow}), Proposition \ref{s-holderbelow} implies
\begin{equation}
|s^+(y_1)-s^+(y_0)|
\ge \cos{\bar{\theta}(y_1)}\hat{C}^-_V|y_1-y_0|^{4n-1}.\nonumber\qquad\square
\end{equation}

{\em Proof of Theorem \ref{holder2}: }
Define $y_i:=t^+(x_i)\in V\subset T_2$. If $y_0\in\Theta(y_1,\theta)$, with $\theta>\bar{\theta}(y_1)$, we have
\begin{displaymath}
|y_1-y_0|<K|x_1-x_0|.
\end{displaymath}
Otherwise, by the uniform continuity of $t^+$, taking $x_0$ sufficiently close to $x_1$, we have $|t^+(x_1)-t^+(x_0)|<\gamma_{\delta_\epsilon}$ and we can apply Proposition \ref{s+holderbelow} to $y_i=t^+(x_i)$, $i=1,0$ to conclude
\begin{displaymath}
|y_1-y_0|<\frac{1}{\hat{C}^+_V}|x_1-x_0|^{\frac{1}{4n-1}}.\qquad\square
\end{displaymath}



\end{document}